\documentclass[11pt]{article}
\usepackage{comment}
\usepackage{graphics,graphicx}
\usepackage{amsfonts,amssymb,amsmath,xcolor,mathrsfs, dsfont,amstext}
\usepackage{amsbsy, amsopn, amscd, amsxtra, amsthm,authblk, enumerate}
\usepackage{upref}
\usepackage{geometry}
\usepackage[displaymath]{lineno}
\usepackage{float}
\usepackage{bbm}
\usepackage{enumitem}
\usepackage{extarrows}
\usepackage{yhmath}
\DeclareSymbolFont{largesymbol}{OMX}{yhex}{m}{n}
\DeclareMathAccent{\Widehat}{\mathord}{largesymbol}{"62}
\numberwithin{equation}{section}

\usepackage[colorlinks,
            linkcolor=blue,
            anchorcolor=green,
            citecolor=blue
            ]{hyperref}

\title{
Uniqueness of generalized conformal restriction measures}
\author{Zhuoyan Xie}
\date{}

\newcommand{\R}{\mathbbm{R}}

\newcommand{\D}{\mathbbm{D}}
\newcommand{\C}{\mathbbm{C}}

\newcommand{\E}{\mathbbm{E}}

\renewcommand{\P}{\mathbbm{P}}

\newcommand{\1}{{\mathbf{1}}}
\newcommand{\hH}{\mathbbm{H}}
\newcommand{\PPP}{\mathbf{P}}
\newcommand{\QQQ}{\mathbf{Q}}

\DeclareMathOperator{\SLE}{SLE}

\newcommand{\wt}{\widetilde}
\newcommand{\balpha}{\boldsymbol{\alpha}}
\newcommand{\bbeta}{\boldsymbol{\beta}}
\newcommand{\baa}{\boldsymbol{a}}
\newcommand{\bxx}{\boldsymbol{x}}
\newcommand{\byy}{\boldsymbol{y}}
\newcommand{\bX}{\boldsymbol{X}}
\newcommand{\bY}{\boldsymbol{Y}}

\newcommand{\bcc}{\boldsymbol{c}}

\newtheorem{theorem}{Theorem}[section]
\newtheorem{definition}[theorem]{Definition}
\newtheorem{lemma}[theorem]{Lemma}

\newtheorem{proposition}[theorem]{Proposition}
\newtheorem{remark}[theorem]{Remark}

\newcommand\wh[1]{\widehat{#1}}

\begin{document}
\maketitle

\begin{abstract}
In this note, we give a new proof of the uniqueness of generalized chordal and radial conformal restriction measures for $c\le 1$, where the $c<-2$ case is unproved before. Our approach (called the SLE arc test) is inspired by the completely conformally invariant (CCI) measure introduced by Lawler and Werner (2000), and also by Cai's recent Green's function test for the Brownian loop catcher (2026). 
\end{abstract}

\section{Introduction}
\label{section-intro}

\subsection{Background}
\label{subsection-background}
In the study of Mandelbrot's conjecture (i.e. the dimension of planar Brownian frontier is $\frac{4}{3}$), Lawler and Werner introduced the completely conformally invariant (CCI) measure which satisfies certain conformal restriction property, and showed that its intersection exponents are the same as the Brownian intersection exponents~\cite{Lawler-Werner-universality}. 
After Schramm's introduction of Schramm--Loewner evolution (SLE)~\cite{schramm0}, Lawler, Schramm and Werner used the $\SLE_6$ exponents and the universality arguments to compute all the Brownian exponents, which led to the proof of Mandelbrot's conjecture~\cite{lsw-bm-exponents1,lsw-bm-exponents2,lsw-bm-exponents3,lsw-analyticity}.

To further understand the relation between planar Brownian motion and $\SLE$,
they~\cite{lsw-restriction} studied the chordal conformal restriction measures and described them in terms of $\SLE_{8/3}$ and its variants. There is also a loop version of conformal restriction measure, constructed by taking the outer boundary of a Brownian loop on the plane~\cite{werner-loops}.
The radial and tri-chordal cases of conformal restriction measures were later studied in~\cite{Wu-radial-restriction,Qian-trichordal}.

For $\kappa\in(0,4]$, $\SLE_\kappa$ satisfies a generalized conformal restriction property involving Brownian loop mass~\cite{lsw-restriction,lawler-werner-soup,dubedat-rho,qian2021generalized}. The generalized conformal restriction measure has a parameter $c\le1$ called central charge, which is related to $\kappa$ by
\begin{equation}
\label{eq-c-kappa}
    c=1-6\left(\frac{\sqrt{\kappa}}{2}-\frac{2}{\sqrt{\kappa}}\right)^2.
\end{equation}
The $c=0$ (or $\kappa=8/3$) case corresponds to the classical conformal restriction measure.
In~\cite{qian2021generalized}, Qian gave the construction of the generalized chordal, radial and trichordal restriction measures for $c\le1$.
The loop counterpart of generalized conformal restriction measure for $c\le1$ is called the Malliavin--Kontsevich--Suhov (MKS) measure~\cite{kontsevich2006}. 
The construction of MKS measure was given by the loop version of $\SLE_\kappa$; see~\cite{BD-SLE2-loop,zhan-loop-measures,werner-sphere-cle}.

Given the existence of the above generalized conformal restriction measures, the uniqueness of generalized conformal restriction measures for $c\le1$ remained open for many years. In~\cite{BJ-cft-sle}, 
Baverez and Jego developed the conformal field theory (CFT) of the SLE loop measure, and provided the first proof of the uniqueness of MKS measure for all $c\le1$. In~\cite{CG25}, Cai and Gao gave a ``probabilistic'' proof of uniqueness for $c\in(0,1]$ by
relating the conformal restriction measures to the Brownian loop soup clusters. Recently, Cai~\cite{cai2026brownianloopcatcher} studied the negative central charge (with $c\in[-2,0)$) extension of Brownian loop soup clusters, which is called the $\lambda$-Brownian loop catcher with $\lambda=-\frac{c}{2}$. To show the uniqueness of $\lambda$-Brownian loop catcher, Cai developed the Green's function test method, which also yields the uniqueness of generalized chordal restriction measure for $c\in[-2,0)$.

In this note, we provide a new proof of the uniqueness of generalized chordal or radial conformal restriction measures for $c\le1$, where $c<-2$ case is unproved before.
Our approach (called \textbf{the SLE arc test}) was inspired by the CCI measure~\cite{Lawler-Werner-universality} and Green's function test method~\cite{cai2026brownianloopcatcher} (see discussions around~\eqref{eq-explain-1} and~\eqref{eq-explain-2}). 
Our method may also have potential applications to give another proof of the uniqueness of MKS measures; see discussions in Section~\ref{section-application}.

\subsection{Main results}

We now introduce the chordal hull and state our main result. Consider the upper half plane $\hH:=\{z\in\C:\Im(z)>0\}$.
Following~\cite{lsw-restriction},
let $\mathcal{H}$ (resp. $\mathcal{H}_+$) be the collection of unbounded closed connected subsets $K$ of $\overline{\hH}$ such that $K\cap\R=\{0\}$ (resp. $K\cap\R=(-\infty,0]$) and $\C\setminus K$ (resp. $\hH\setminus K$) is connected. We call an $H\in\mathcal{H}$ (resp. $H\in\mathcal{H}_+$) a two-sided hull (resp. one-sided hull). Let $\mathcal{Q}$ be the collection of compact sets $A\subset\overline{\hH}$ such that $\hH\setminus A$ is simply connected and $0\notin A$, and let $\mathcal{Q}_+$ (resp. $\mathcal{Q}_-$) be the collection of $A\in\mathcal{Q}$ such that $A\cap(-\infty,0]=\emptyset$ (resp. $A\cap[0,+\infty)=\emptyset$). For $\mathcal{H}$ (resp. $\mathcal{H}_+$), we equip it with the $\sigma$-algebra generated by $\{H:H\cap A=\emptyset\}$ for all $A\in\mathcal{Q}$ (resp. $A\in\mathcal{Q}_+$).

For two closed subsets $A_1,A_2\subset\overline{\hH}$, let $m_{\hH}(A_1,A_2)$ be the total mass of the Brownian loops on $\hH$ that intersect both $A_1$ and $A_2$ under the Brownian loop measure $\mu^{\text{BL}}$, which is defined in~\cite{lawler-werner-soup}. The $c\in[-2,0]$ case of the following result was first proved in~\cite[Theorem 1.6]{cai2026brownianloopcatcher}. We now extend its range to $c\le1$.
\begin{theorem}
    \label{thm-chordal-hull}
    Fix $c\le 1$. For a probability law $\PPP$, we use $\E_\PPP$ to denote the expectation under $\PPP$.\\
    (a) Suppose $\PPP_1,\PPP_2$ are probability measures on $H\in\mathcal{H}_+$ satisfying
    \begin{equation*}
\E_{\PPP_1}[\1_{\{H\cap A=\emptyset\}}e^{\frac{c}{2}m_{\hH}(H,A)}]=\E_{\PPP_2}[\1_{\{H\cap A=\emptyset\}}e^{\frac{c}{2}m_{\hH}(H,A)}]\in[0,+\infty]
\end{equation*}
    for every $A\in\mathcal{Q}_+$. Then $\PPP_1=\PPP_2$. \\
    (b) Suppose $\PPP_1,\PPP_2$ are probability measures on $H\in\mathcal{H}$ satisfying
    \begin{equation*}
\E_{\PPP_1}[\1_{\{H\cap A=\emptyset\}}e^{\frac{c}{2}m_{\hH}(H,A)}]=\E_{\PPP_2}[\1_{\{H\cap A=\emptyset\}}e^{\frac{c}{2}m_{\hH}(H,A)}]\in[0,+\infty]
\end{equation*}
    for every $A\in\mathcal{Q}$. Then $\PPP_1=\PPP_2$.
\end{theorem}

For $c<0$, Theorem~\ref{thm-chordal-hull} means that the law of the hull $H$ is determined by the law of hull obtained by adding all the Brownian loops that intersect $H$ from an independent Brownian loop soup with intensity $-\frac{c}{2}\mu^{\text{BL}}$. In particular, Theorem~\ref{thm-chordal-hull} solves Question 2 in~\cite[Section 8.3]{werner-restriction-survey}. 
\begin{theorem}
    \label{cor-uniquness-loop-decorator}
    Suppose $H_1,H_2$ are two random one-sided (or two-sided) hulls on $\hH$. For $c>0$, let $\mathcal{L}_\hH^{\frac{c}{2}}$ be another independent Brownian loop soup on $\hH$ with intensity measure $\frac{c}{2}\mu^{\text{BL}}$. Let $L(H_1)$ and $L(H_2)$ be the one-sided (or two-sided) hulls obtained by filling all the bounded complementary components of 
    \begin{equation*}
\overline{H_1\cup\bigcup_{\{\ell\in\mathcal{L}_\hH^{c/2},\ell\cap H_1\neq\emptyset\}}\ell}\quad\text{ and }\quad\overline{H_2\cup\bigcup_{\{\ell\in\mathcal{L}_\hH^{c/2},\ell\cap H_2\neq\emptyset\}}\ell}
\end{equation*}
    respectively. If $L(H_1)\overset{d}{=}L(H_2)$, then $H_1\overset{d}{=}H_2$.
\end{theorem}
\begin{proof}
    We denote the law of $H_1$ (resp. $H_2$) by $\PPP_1$ (resp. $\PPP_2$) and the law of $\mathcal{L}_\hH^{\frac{c}{2}}$ by $\P$. For any $A\in\mathcal{Q_+}$ (if $H_1,H_2$ are one-sided hulls) or $A\in\mathcal{Q}$ (if $H_1,H_2$ are two-sided hulls), we have for $i=1,2$,
    \begin{equation*}
(\P\times\PPP_i)[L(H_i)\cap A=\emptyset]=\E_{\PPP_i}[\1_{\{H_i\cap A=\emptyset\}}e^{-\frac{c}{2}m_{\hH}(H_i,A)}].
\end{equation*}
Therefore, $\E_{\PPP_1}[\1_{\{H_1\cap A=\emptyset\}}e^{-\frac{c}{2}m_{\hH}(H_1,A)}]=\E_{\PPP_2}[\1_{\{H_2\cap A=\emptyset\}}e^{-\frac{c}{2}m_{\hH}(H_2,A)}]$. By Theorem~\ref{thm-chordal-hull}, we obtain that $\PPP_1=\PPP_2$.
\end{proof}

Note that for $c\le0$, the number $\E_{\PPP_1}[\1_{\{H\cap A=\emptyset\}}e^{\frac{c}{2}m_{\hH}(H,A)}]\in[0,1]$, while for $c\in(0,1]$ it may be infinity. The $c=0$ case of Theorem~\ref{thm-chordal-hull} follows from the definition of the $\sigma$-algebra on $\mathcal{H}$ or $\mathcal{H}_+$.
For $c\in(0,1]$, we can relate $\E_{\PPP_1}[\1_{\{H\cap A=\emptyset\}}e^{\frac{c}{2}m_{\hH}(H,A)}]$ to the Brownian loop soup clusters with intensity $\frac{c}{2}\mu^{\text{BL}}$, as in~\cite{CG25}\footnote{They prove the uniqueness of generalized chordal restriction measures for $c\in(0,1]$, which is weaker than the $c\in(0,1]$ case of Theorem~\ref{thm-chordal-hull}. Nevertheless, their method can be applied to the $c\in(0,1]$ case of Theorem~\ref{thm-chordal-hull}.}. The $c\in[-2,0)$ case was proved in~\cite[Theorem 1.6]{cai2026brownianloopcatcher} by the Green's function test method, which uses $\lambda$-Brownian loop-catcher ($\lambda=-\frac{c}{2}$) for multiple Brownian paths to show that the joint moments of Poisson kernels of $\hH\setminus H$ agree under $\PPP_1$ or $\PPP_2$, and hence $\PPP_1=\PPP_2$; see~\cite[Section 4]{cai2026brownianloopcatcher}. Note that this method does not cover $c<-2$, since there is no $\lambda$-Brownian loop catcher for $\lambda=-\frac{c}{2}>1$ by~\cite[Theorem 1.1]{cai2026brownianloopcatcher}.

One feature of Green's function test method is to take the expectation of the random test sets, which is similar to taking $F=1$ in~\eqref{eq-explain-1}, but with different random test sets there. 
It is natural to ask, how much information is lost if we directly take $F=1$ in~\eqref{eq-explain-1}?
Indeed, if we allow the arbitrariness of $F$, then~\eqref{eq-explain-1} is equivalent to the condition in Theorem~\ref{thm-chordal-hull}, and hence no information is lost.
Our approach here uses the arbitrariness of $F$, which can be applied to all the $c\le1$. 
Inspired by the CCI measure~\cite{Lawler-Werner-universality},
we define a well-chosen SLE arc measure on $\hH$ (say $\nu$) such that this arc separates 1 from $\infty$ with two endpoints lying on $\R_+$, see~\eqref{eq-SLE-arc}. 
We take $A$ to be the closure of bounded component of $\hH\setminus\eta$, where $\eta$ is an SLE arc sampled from $\nu$. Therefore, in the setting of Theorem~\ref{thm-chordal-hull} (a), for any non-negative measurable function $F$ on $\eta$, we have~\footnote{The following argument is a simplified version. The measure $\nu$ here is an infinite measure and needs to be normalized as a probability measure. The rigorous argument is given in the proof of Theorem~\ref{thm-chordal-hull} (a).}
\begin{equation}
\label{eq-explain-1}
(\nu\times \E_{\PPP_1})[\1_{\{H\cap \eta=\emptyset\}}e^{\frac{c}{2}m_\hH(H,\eta)}F(\eta)]=(\nu\times \E_{\PPP_2})[\1_{\{H\cap \eta=\emptyset\}}e^{\frac{c}{2}m_\hH(H,\eta)}F(\eta)].
\end{equation}

As we will see, this SLE arc measure $\nu$ satisfies some generalized conformal restriction property (Lemma~\ref{lem-restriction-SLE-arc}). Therefore, using Tonelli's theorem on~\eqref{eq-explain-1},
\begin{equation}
\label{eq-explain-2}
\int\left(\int F(\eta)\nu_{\hH\setminus H}(d\eta)\right)\PPP_1(dH)=\int\left(\int F(\eta)\nu_{\hH\setminus H}(d\eta)\right)\PPP_2(dH).
\end{equation}
Here $\nu_{\hH\setminus H}$ is the pushforward of $\nu$ under some conformal map from $\hH$ to $\hH\setminus H$. Note that~\eqref{eq-explain-2} implies that the marginal law of $\eta$ is the same under $\nu_{\hH\setminus H}(d\eta)\PPP_1(dH)$ or $\nu_{\hH\setminus H}(d\eta)\PPP_2(dH)$. To convert this fact into a more precise statement, we can use conformal maps to parametrize domains.
For $H_1$ (resp. $H_2$) sampled from $\PPP_1$ (resp. $\PPP_2$), let $f_1$ (resp. $f_2$) be the conformal map from $\hH$ to $\hH\setminus H_1$ (resp. $\hH\setminus H_2$) with $0,1,\infty$ fixed; let $g$ be the conformal map from $\hH$ to the bounded component of $\hH\setminus\eta$ with 1 fixed, where $\eta$ is sampled from $\nu$ (independent of $\PPP_1,\PPP_2$). Thanks to the Schwarz reflection theorem, these conformal maps can be extended to univalent maps on $\C\setminus(-\infty,0]$ (by doubling the domain along $\R_+$).
In particular, ~\eqref{eq-explain-2} implies that $f_1\circ g\overset{d}{=}f_2\circ g$. We then show in Section~\ref{section-composition} that this identity actually forces $f_1\overset{d}{=}f_2$.

A direct application of Theorem~\ref{thm-chordal-hull} is the uniqueness of generalized chordal conformal restriction measure, which we now recall; see~\cite{lsw-restriction} and~\cite[Definition 6.1]{qian2021generalized}. 
\begin{definition}[generalized chordal restriction measure]
    \label{def-chordal-restriction}
    Fix $c\le 1$ and $\alpha\in\R$. A probability measure $\PPP$ on $\mathcal{H}$ (resp. $\mathcal{H}_+$) is called a generalized two-sided (resp. one-sided) chordal restriction measure with central charge $c$ and exponent $\alpha$, if it satisfies
    \begin{equation*}
    \1_{\{H\cap A=\emptyset\}}\PPP(dH)=(\phi_A'(0)\phi_A'(\infty))^\alpha e^{-\frac{c}{2}m_{\hH}(H,A)}\PPP_A(dH)
\end{equation*}
for all $A\in\mathcal{Q}$ (resp. $A\in\mathcal{Q}_+$). Here $\phi_A:\hH\setminus A\to\hH$ is conformal with $0,\infty$ fixed, and $\phi_A'(\infty):=\lim\limits_{z\to\infty}\frac{z}{\phi_A(z)}$.\footnote{This limit exists and belongs to $(0,\infty)$, see for example~\cite[Section 2]{lsw-restriction}.} The value $\phi_A'(0)\phi_A'(\infty)\in(0,1]$\footnote{This follows from~\cite[Section 2, equation (2.4)]{lsw-restriction}.} does not rely on the choice of $\phi_A$. The law $\PPP_A$ is defined as the pushforward of $\PPP$ under $\phi_A^{-1}$.
\end{definition}
\begin{theorem}
    \label{thm-unique-chordal-restriction}
    Fix $c\le 1$ and $\kappa\in(0,4]$ related to $c$ via~\eqref{eq-c-kappa}.
    For each $\alpha\in\R$, there is at most one generalized two-sided (resp. one-sided) chordal restriction measure with central charge $c$ and exponent $\alpha$. Moreover, the generalized two-sided (resp. one-sided) chordal restriction measure with central charge $c$ exists if and only if $\alpha\ge\frac{6-\kappa}{2\kappa}$ (resp. $\alpha\ge0$). 
\end{theorem}

The $c=0$ case of Theorem~\ref{thm-unique-chordal-restriction} case was studied in~\cite{lsw-restriction}, the $c\in(0,1]$ case was treated in~\cite{CG25} and the $c\in[-2,0)$ case was shown in~\cite{cai2026brownianloopcatcher}. 
Also\footnote{This was suggested by Gefei Cai in private communications.}, the radial probes in~\cite[Section 4.3]{cai2026brownianloopcatcher} can give another proof of Theorem~\ref{thm-unique-chordal-restriction} for $c<0$.
The construction of generalized chordal restriction measure for all $c\le 1$ was given in~\cite[Proposition 6.2]{qian2021generalized}.
Note that our SLE arc test method can also be applied to other topologies, such as the uniqueness of radial hull (Theorem~\ref{thm-radial-hull}) or tri-chordal hull (Remark~\ref{remark-trichordal}), see Section~\ref{section-arc} for more details.

\subsection{Discussions, applications and outlooks}
\label{section-application}

As pointed out in~\cite[Section 8]{werner-restriction-survey},
the uniqueness of conformal restriction measure has applications to the reversibility and duality of SLE. 
For the SLE loop measure, its reversibility and duality was shown in~\cite[Section 7]{BJ-cft-sle}; see also~\cite{ACSW24a} for another proof of the duality of SLE loop measure.

Here, we use Theorem~\ref{thm-unique-chordal-restriction} to give another proof of the reversibility of chordal SLE. We briefly review the history of the reversibility of chordal SLE. It was first conjectured in~\cite{schramm-sle} that chordal $\SLE_\kappa$ is reversible for $0<\kappa\le 8$. After this, Zhan proved reversibility for $0<\kappa\le 4$~\cite{zhan-reversibility}, and Miller and Sheffield proved the range $4<\kappa<8$ by the imaginary geometry~\cite{ig3}. The Corollary~\ref{cor-reverse-SLE} gives another proof of~\cite{zhan-reversibility}. The imaginary geometry approach can be used to study the reversibility of $\SLE$ with more forces points~\cite{ig2} and non-simple case $\kappa>4$~\cite{ig3,ig4}, while the conformal restriction approach cannot cover these ranges.
\begin{theorem}
    \label{cor-reverse-SLE}
    Fix $c\le 1$ and $\kappa\in(0,4]$ related to $c$ via~\eqref{eq-c-kappa}. Suppose $\PPP$ is a generalized one-sided conformal restriction measure with exponent $\alpha>0$ and central charge $c$. Let $I(z):=1/{\overline{z}}$ be the reflection on $\hH$. Then the pushforward measure $I_\#\PPP$ is still a generalized one-sided conformal restriction measure with exponent $\alpha>0$ and central charge $c$. In particular, we have $\PPP=I_\#\PPP$.

As a corollary, for an $\SLE_\kappa(\rho)$ ($\rho>-2$) curve $\eta$ on $\hH$ from 0 to $\infty$ with force point $0_-$, we have $\eta\overset{d}{=}1/\overline{\eta}$, where we do a time-reverse and also a time-change to $1/\overline{\eta}$.
\end{theorem}
\begin{proof}
Note that $\mathcal{Q}_+$ is invariant under $I$. For any $A\in\mathcal{Q}_+$,
we have $\phi_A'(0)\phi_A'(\infty)=\phi_{I(A)}'(0)\phi_{I(A)}'(\infty)$. By the inversion invariance of Brownian loop measure on $\hH$, we have $m_\hH(H,A)=m_\hH(I(H),I(A))$ for $H\in\mathcal{H}_+$ and $A\in\mathcal{Q}_+$ with $H\cap A=\emptyset$. Combined, we can verify that $I_\#\PPP$ satisfies the generalized one-sided conformal restriction measure with exponent $\alpha>0$ and central charge $c$. By Theorem~\ref{thm-unique-chordal-restriction}, we have $\PPP=I_\#\PPP$. Finally, $\SLE_\kappa(\rho)$ with $\rho>-2$ from 0 to $\infty$ with force point $0_-$ gives the one-sided restriction measure with exponent $\alpha=\frac{(\rho+2)(\rho+6-\kappa)}{4\kappa}$ for all $\alpha>0$; see~\cite{lsw-restriction,dubedat-rho} and also~\cite[Section 6.2.1]{qian2021generalized}.
\end{proof}

Recently, Baverez and Jego~\cite{BJ-cft-sle} provided the first proof of the
uniqueness of MKS measure (up to a multiplicative constant) for all $c\le1$, which is an application of their development
of the conformal field theory (CFT) of the SLE loop measure.
Their work goes beyond the uniqueness question, including the
representations of Virasoro algebras, a geometric interpretation
of the Shapovalov form via the integration-by-parts formula, the spectral resolution of this Shapovalov form, and therefore revealing a rich
algebraic structure underlying SLE loop measures. See also the independent and contemporaneous work~\cite{GQW-conformal-restriction} on the natural Virasoro representations from the SLE loop measure.

As mentioned in Section~\ref{subsection-background}, our method 
may also lead to another proof of the uniqueness of the MKS measure.  
The idea is to test an
unknown MKS measure against the SLE loop measure. To obtain
finite measures, one possible way is to restrict the SLE loop measure to
loops whose conformal radii lie in a fixed interval, and reweight the MKS measure by the conformal radius,
in the spirit of Lemma~\ref{lem-finite-SLE-arc}. However, further arguments are needed to make this idea rigorous.

We collect here some other possible applications and outlooks. 
\begin{itemize}
    \item By testing against multiple SLE loops, we can
apply Proposition~\ref{prop-joint-composition} to extend the
range of~\cite[Theorem 1.5]{cai2026brownianloopcatcher}
to all $\lambda>0$. Note that our proof of the uniqueness of chordal two-sided hull has already used this multiple flavor; see also Remark~\ref{remark-trichordal}.
This extension is also related to the Question 4
in~\cite[Section 1.6]{cai2026brownianloopcatcher}.
We plan to address these questions in the future.
\item For $c\in(0,1]$ case, Cai-Gao's proof gives the non-intersection probability of generalized chordal conformal restriction measure in terms of the Brownian loop soup clusters; see e.g.~\cite[Theorem 1.3]{CG25}.
It is then interesting to find some probabilistic reversal formulas for $c<0$.
\item The Green's function test can be applied to the three-dimensional case; see the recent work on the characterization of the scaling limit of three-dimensional loop-erased random walk~\cite{Cai-3D-LERW}. It would be interesting to see whether the strategy here has an analogue in higher dimensions.
\end{itemize}

\noindent{\bf Acknowledgment.} 
I'd like to thank Gefei Cai for explaining his new work \textit{The Brownian loop catcher} to me in July 2026, especially the Green's function test. I'd like to thank Gefei Cai, Shengzhi Jin, Xin Sun and Baojun Wu for many helpful discussions.
This work was supported by National Key R\&D Program of China (No.\ 2023YFA1010700) and National Natural Science Foundation of China (Grant No. 12526204).\\
~\\
\noindent{\bf AI usage.} I developed the proof strategy, including the SLE arc measure and SLE arc test method demonstrated in~\eqref{eq-explain-1},~\eqref{eq-explain-2} and Section~\ref{section-arc}.
Within this strategy, I used AI-tool to explore and find a short proof of Lemma~\ref{lem-linear-algebra-uniqueness}. I also used AI-tool to assist to give proofs in Appendix~\ref{appendix-measurable}. I take full responsibility for the correctness of the proof. In the final stage, I used AI-tool to polish the language.

\section{Random composition: a common toolbox}
\label{section-composition}
Let $\D:=\{z\in\C:|z|<1\}$ be the unit disk.
Let $\mathcal{U}$ be the space of univalent map (i.e. injective and holomorphic) $f:\D\to\C$ with $f(0)=0$ and $f'(0)>0$. For $M>0$, let $\mathcal{U}_M:=\{f\in\mathcal{U}:|f'(0)|\le M\}$. Let $\mathcal{D}:=\{f\in\mathcal{U}:f(\D)\subset\D\}$. 
Note that Schwarz's Lemma implies that $\mathcal{D}\subset\mathcal{U}_1$.
For $\mathcal{U},\mathcal{U}_M,\mathcal{D}$, we equip them with the Borel $\sigma$-algebra generated by the topology of locally uniform convergence on $\D$.
\begin{proposition}
    \label{prop-unique-composition}
    Fix a probability measure $\QQQ$ on $\mathcal{D}$ and $M>0$.
    Suppose $\PPP,\wt\PPP$ are two probability measures on $\mathcal{U}_M$.
    Let $(f,g)$ be sampled from $\PPP\times\QQQ$ and $(\wt f,\wt g)$ be sampled from $\wt\PPP\times\QQQ$. 
    If $f\circ g\overset{d}{=}\wt f\circ\wt g$, then $\PPP=\wt\PPP$.
\end{proposition}

Note that if $g$ is deterministic (or $\QQQ=\delta_g$), then $f\circ g{=}\wt f\circ g$ implies that $f=\wt f$, since $g(\D)$ is an open set. We also need a multiple version of Proposition~\ref{prop-unique-composition}. For an integer $n\ge1$, let $\mathcal{U}_M^{\otimes n}:=\mathcal{U}_M\times\cdots\times\mathcal{U}_M$ be the $n$-fold product space of $\mathcal{U}_M$, and $\QQQ^{\otimes n}:=\QQQ\times\cdots\times\QQQ$ be the $n$-fold product measure of $\QQQ$.
\begin{proposition}
    \label{prop-joint-composition}
    Fix an integer $n\ge1$, a number $M>0$ and a probability measure $\QQQ$ on $\mathcal{D}$.
    Suppose $\PPP$ and $\wt\PPP$ are two probability measures on $\mathcal{U}_M^{\otimes n}$. Let $((f_1,\cdots,f_n),(g_1,\cdots,g_n))$ be sampled from $\PPP\times\QQQ^{\otimes n}$ and $((\wt f_1,\cdots,\wt f_n),(\wt g_1,\cdots,\wt g_n))$ be sampled from $\wt\PPP\times\QQQ^{\otimes n}$. If
    \begin{equation*}
\left(f_1\circ g_1,\cdots,f_n\circ g_n\right)\overset{d}{=}\left(\wt f_1\circ \wt g_1,\cdots,\wt f_n\circ \wt g_n\right),\text{ then }\PPP=\wt\PPP.
\end{equation*}
\end{proposition}

To show Proposition~\ref{prop-unique-composition}, we look at the coefficients of the conformal maps $f,g$ and $f\circ g$, and show that the joint moments of the coefficients in $f$ is determined. The idea of \textit{using the moments of the bounded coefficients to determine the law the random conformal map} was previously appeared in the proof of the uniqueness of MKS measure~\cite{BJ-cft-sle}; see also~\cite{Werner-measure-uniqueness}. Our proof of Proposition~\ref{prop-unique-composition} also used this idea, but with different mechanics on how to determine the joint moments of coefficients.
More precisely, if $f(z)=\sum_{n\ge1}a_nz^n,g(z)=\sum_{n\ge1}b_nz^n,f(g(z))=\sum_{n\ge1}c_nz^n$, then there is a $n\times n$ lower triangular matrix $B$ whose elements only depends on $\{b_n\}_{n\ge1}$, such that $\bcc=B\baa$ where $\bcc:=(c_1,\cdots,c_n)^T,\baa:=(a_1,\cdots,a_n)^T$. We have the following result.
\begin{lemma}
    \label{lem-linear-algebra-uniqueness}
    Fix integer $n\ge1$. Suppose $B=(b_{ij})_{1\le i,j\le n}\in\C^{n\times n}$ is a random $n\times n$ complex matrix, such that $B$ is bounded (i.e. there is $M>0$ such that $|b_{ij}|<M$ holds for all $i,j$ with probability 1), lower triangular (i.e. $b_{ij}=0$ if $i<j$), and $b_{ii}>0$ for all $i$. Suppose $\bxx,\byy\in\C^n$ are two bounded random complex vectors independent of $B$. If $B\bxx\overset{d}{=}B\byy$, then $\bxx\overset{d}{=}\byy$.
\end{lemma}

We now prove these results.
\begin{proof}[Proof of Lemma~\ref{lem-linear-algebra-uniqueness}]
    Since $\bxx,\byy$ are bounded, it suffices to verify that $\E[\bxx^{\balpha}\overline{\bxx}^{\bbeta}]=\E[\byy^{\balpha}\overline{\byy}^{\bbeta}]$ holds for all $\balpha,\bbeta\in\mathbbm{N}^n$. Here, we use $\bxx^{\balpha}:=x_1^{\alpha_1}\cdots x_n^{\alpha_n},\overline{\bxx}^{\bbeta}:=\overline{x_1}^{\beta_1}\cdots \overline{x_n}^{\beta_n}$, and similarly define $\byy^{\balpha},\overline\byy^{\bbeta}$. 
    For a monomial $\bxx^{\balpha}\overline{\bxx}^{\bbeta}$, we define $d:=\sum_{k=1}^n k(\alpha_k+\beta_k)$ to be its index.
    We prove this by induction on the index $d$.
    We write $\bX:=B\bxx,\bY:=B\byy$. For $d=0$, both expectations are 1.
    For $d=1$, $\E[b_{11}]\E[x_1]=\E[X_1]=\E[Y_1]=\E[b_{11}]\E[y_1]$ and $\E[b_{11}]>0$ implies that $\E[x_1]=\E[y_1]$. For general $d$, we have the expansion
    \begin{equation*}
    \begin{aligned}
&\E\left[\bX^{\balpha}\overline{\bX}^{\bbeta}\right]=\E\left[\prod_{k=1}^n\left(\sum_{i=1}^kb_{ki}x_i\right)^{\alpha_k}\left(\sum_{j=1}^k\overline{b_{kj}x_j}\right)^{\beta_k}\right]\\
&=\E\left[\prod_{k=1}^nb_{kk}^{\alpha_k+\beta_k}\right]\E\left[\bxx^{\balpha}\overline{\bxx}^{\bbeta}\right]+\sum_{(\balpha',\bbeta'):\text{ index }\bxx^{\balpha'}\overline{\bxx}^{\bbeta'}<d}C_{\balpha',\bbeta'}\E[\bxx^{\balpha'}\overline{\bxx}^{\bbeta'}].
\end{aligned}
\end{equation*}
By induction,
$\E\left[\bX^{\balpha}\overline{\bX}^{\bbeta}\right]=\E\left[\bY^{\balpha}\overline{\bY}^{\bbeta}\right]$, $\E\left[\prod_{k=1}^nb_{kk}^{\alpha_k+\beta_k}\right]>0$ implies that $\E[\bxx^{\balpha}\overline{\bxx}^{\bbeta}]=\E[\byy^{\balpha}\overline{\byy}^{\bbeta}]$.
\end{proof}

\begin{proof}[Proof of Proposition~\ref{prop-unique-composition}]

    For $f\in\mathcal{U}_M$ and $g\in\mathcal{D}$, we have the expansion $f(z)=\sum_{n\ge1}a_nz^n$ and $g(z)=\sum_{n\ge1}b_nz^n$ with $a_1,b_1>0$. We have $a_1\le M$, $b_1\le 1$ (by Schwarz's lemma) and $|a_n|\le Mn,|b_n|\le n$ for $n\ge2$ by de Branges~\cite{deBranges1985}. Let $f(g(z))=\sum_{n\ge1}c_nz^n$. For each $n\ge1$:
\begin{equation}
\label{eq-proof-appendix-1}
\begin{pmatrix}
c_1\\ \vdots\\ c_n
\end{pmatrix}
=
\begin{pmatrix}
b_1 & & 0\\
& \ddots &\\
\ast & & b_1^n
\end{pmatrix}
\begin{pmatrix}
a_1\\ \vdots\\ a_n
\end{pmatrix}.
\end{equation}
Here the matrix in~\eqref{eq-proof-appendix-1} is lower triangular with positive diagonals, and $\ast$ are polynomials composed by $(b_k)_{1\le k\le n}$, which are bounded. Furthermore, this matrix only depends on $g$ and are independent of $f$.
By Lemma~\ref{lem-linear-algebra-uniqueness}, the law of $(a_1,\cdots,a_n)$ is uniquely determined for each $n\ge1$, and hence the law of $f$ is uniquely determined. Therefore, $\PPP=\wt\PPP$.
\end{proof}

\begin{proof}[Proof of Proposition~\ref{prop-joint-composition}]
For $1\le j\le n$, we have the expansions $f_j(z)=\sum_{k\ge1}a_{j,k}z^k,f_j(g_j(z))=\sum_{k\ge1}c_{j,k}z^k$.
For each $N\ge1$, let
$\baa_j:=(a_{j,1},\ldots,a_{j,N})^T$ and
$\bcc_j:=(c_{j,1},\ldots,c_{j,N})^T$.
As in the proof of Proposition~\ref{prop-unique-composition}, we have
\begin{equation*}
\begin{pmatrix}
\bcc_1\\ \vdots\\ \bcc_n
\end{pmatrix}
=\text{diag}(B_{1,N},\cdots,B_{n,N})
\begin{pmatrix}
\baa_1\\ \vdots\\ \baa_n
\end{pmatrix},
\end{equation*}
where $B_{j,N}$ is the composition matrix associated with $g_j$ as in~\eqref{eq-proof-appendix-1}.
Here $\text{diag}(B_{1,N},\cdots,B_{n,N})$ is bounded and lower triangular
with positive diagonals. Furthermore, it only depends on
$(g_1,\ldots,g_n)$ and is independent of $(f_1,\ldots,f_n)$.
By Lemma~\ref{lem-linear-algebra-uniqueness}, the joint law of
$(\baa_1,\ldots,\baa_n)$ is uniquely determined for
each $N\ge1$, and hence the joint law of $(f_1,\ldots,f_n)$
is uniquely determined. Therefore, $\PPP=\wt\PPP$.
\end{proof}

As we will see in Section~\ref{section-arc},
the proof of uniqueness of one-sided hull or radial hull uses Proposition~\ref{prop-unique-composition}; while the proof of two-sided hull uses the $n=2$ case of Proposition~\ref{prop-joint-composition}.

\section{SLE arc test}
\label{section-arc}
In this section, we use the SLE arc test method to prove the uniqueness results for chordal hull and radial hull. We assume $\kappa\in(0,4]$ and $c\le 1$ are related via~\eqref{eq-c-kappa}
throughout this section.

\subsection{An SLE arc measure}
In this section, we introduce an SLE arc measure on $\hH$ with marked points $0,1,\infty$ and study its conformal restriction property. For $\kappa\in(0,4]$ and $\rho>-2$,
the $\SLE_\kappa(\rho)$ measure is a variant of $\SLE_\kappa$ measure with a force point, which was first introduced in~\cite{lsw-restriction} to study the one-sided conformal restriction measure; see also~\cite{dubedat-rho}.
For $0<x<y<\infty$, let $\SLE_{x\to y}$ be the $\SLE_{\kappa}(\kappa-2)$ probability measure on $\hH$ from $x$ to $y$ with force point $x_+$. We define the following infinite measure
\begin{equation}
\label{eq-SLE-arc}
\nu_\hH(d\eta):=\iint_{0<x<1<y<\infty}\SLE_{x\to y}(d\eta)\frac{dxdy}{|x-y|^2}.
\end{equation}

We briefly explain the reason that we choose $\rho=\kappa-2$ and the integral kernel $\frac{dxdy}{|x-y|^2}$. For $\kappa=8/3$, this measure is exactly the left boundary of a Brownian excursion measure on $\hH$ with one endpoint in $(0,1)$ and another endpoint in $(1,\infty)$, which is an example of CCI measure~\cite{Lawler-Werner-universality}. The kernel $\frac{dxdy}{|x-y|^2}$ is the boundary Poisson kernel on $\hH$. We need to define a probability measure from the infinite measure $\nu_\hH$. A natural way is to weight this measure by the boundary conformal derivative. As we will explain in Remark~\ref{remark-weight-conformal-radius}, this is also a natural way to change the exponent of a generalized restriction measure.

For $\nu_\hH$-almost every $\eta$, it is a continuous simple curve that intersects $\R$ only at its endpoints.\footnote{Since $\kappa-2>\frac{\kappa}{2}-2$, $\eta$ intersects $\R$ only at its endpoints~\cite[Lemma 8.3]{lsw-restriction}. Moreover,~\cite[Remark 2.3]{ig1} yields that $\eta$ is a continuous simple curve.} Therefore, we can define $D_\eta$ to be the unique bounded connected component of $\hH\setminus\eta$. Note that $\overline{D_\eta}\in\mathcal{Q}_+$. 
Let $g_\eta:\hH\to D_\eta$ be the unique conformal map with $g_\eta(1)=1,g_\eta(0)=\eta_x,g_\eta(\infty)=\eta_y$ where $\eta_x$ (resp. $\eta_y$) is the starting (resp. ending) point of $\eta$. We now introduce a probability measure $\wt{\nu}_{\hH}$ obtained by weighting $\nu_\hH$ by the boundary conformal derivative.
\begin{lemma}
    \label{lem-finite-SLE-arc}
    We have $Z:=\int |g_\eta'(1)|\nu_\hH(d\eta)\in(0,\infty)$. We define $\wt{\nu}_\hH(d\eta):=\frac{1}{Z}|g_\eta'(1)|\nu_\hH(d\eta)$.
\end{lemma}
\begin{proof}
    For $\eta$ with endpoints $0<x<1<y$, we claim that $0<g_\eta'(1)\le1-x$. Indeed, since $\partial D_\eta=\eta\cup[x,y]$ is a Jordan loop, the map $g_\eta$ extends continuously to $\R$, and $g_\eta(\R_+)\subset[x,y]$.
    In particular, we have $g_\eta'(1)>0$.
    Let $\wh{g}_\eta(z):=g_\eta(z)$ if $\Im(z)>0$ and $\wh{g}_\eta(z):=\overline{g_\eta(\overline{z})}$ if $\Im(z)<0$. By Schwarz reflection theorem, $\wh{g}_\eta$ is a univalent map on $\C\setminus(-\infty,0]$ with 1 fixed. 
    Let $J(z)=(\frac{1+z}{1-z})^2$ be a conformal map from $\D$ to $\C\setminus(-\infty,0]$ with $J(0)=1$. By applying Schwarz's Lemma to $J^{-1}\circ(\frac{\wh{g}_\eta(\cdot)-x}{1-x})\circ J$ (which is holomorphic from $\D$ to $\D$ with 0 fixed), we conclude that $0<g_\eta'(1)\le 1-x$.
    
    Therefore, we have
    \begin{equation*}
0<Z\le\iint_{0<x<1<y<\infty}(1-x)\frac{dxdy}{|x-y|^2}=1.\qedhere
\end{equation*}
\end{proof}

The measure $\wt{\nu}_\hH$ satisfies a simple conformal restriction property, which follows from the conformal restriction property of $\SLE_\kappa(\rho)$~\cite{lsw-restriction,dubedat-rho}, see also~\cite[Lemma 4.3]{Wu-CMP-hSLE} for the version we use. For a one-sided hull $H\in\mathcal{H}_+$, let $f_H:\hH\to\hH\setminus H$ be the unique conformal map with $0,1,\infty$ fixed. We define $\wt{\nu}_{\hH\setminus H}$ to be the pushforward of $\wt{\nu}_\hH$ under $f_H$. 
\begin{lemma}
    \label{lem-restriction-SLE-arc}
    Fix a one-sided hull $H\in\mathcal{H}_+$. The following identity holds
    \begin{equation}
    \label{eq-restriction-SLE-arc}
|f_H'(1)|\wt{\nu}_{\hH\setminus H}(d\eta)=\1_{\{\eta\cap H=\emptyset\}}e^{\frac{c}{2}m_\hH(H,\eta)}\wt{\nu}_\hH(d\eta).
\end{equation}
\end{lemma}
\begin{proof}
Recall $\nu_\hH$ in~\eqref{eq-SLE-arc}. We define $\nu_{\hH\setminus H}$ to be the pushforward of $\nu_\hH$ under $f_H$. We first show that
\begin{equation}
\label{eq-restriction-SLE-arc-1}
{\nu}_{\hH\setminus H}(d\eta)=\1_{\{\eta\cap H=\emptyset\}}e^{\frac{c}{2}m_\hH(H,\eta)}{\nu}_\hH(d\eta).
\end{equation}

To show~\eqref{eq-restriction-SLE-arc-1}, we need to prove that for $0<x<1<y<\infty$,
\begin{equation}
\label{eq-restriction-SLE-arc-2}
(f_H)_{\#}\SLE_{x\to y}(d\eta)=\frac{|f'_H(x)f'_H(y)|(x-y)^2}{|f_H(x)-f_H(y)|^2}\1_{\{\eta\cap H=\emptyset\}}e^{\frac{c}{2}m_\hH(H,\eta)}\SLE_{f_H(x)\to f_H(y)}(d\eta).
\end{equation}

Indeed, let $f_{x,y,H}$ be the conformal map from $\hH$ to $\hH\setminus H$ with $x,1,y$ fixed. By~\cite[Lemma 4.3]{Wu-CMP-hSLE}, we have
\begin{equation}
\label{eq-restriction-SLE-arc-3}
(f_{x,y,H})_{\#}\SLE_{x\to y}(d\eta)=|f_{x,y,H}'(x)f_{x,y,H}'(y)|\1_{\{\eta\cap H=\emptyset\}}e^{\frac{c}{2}m_\hH(H,\eta)}\SLE_{x\to y}(d\eta).
\end{equation}
Let $\varphi:\hH\to\hH$ be conformal with $\varphi(f_H(x))=x,\varphi(1)=1,\varphi(f_H(y))=y$. We have $\varphi'(f_H(x))\varphi'(f_H(y))\\
=\frac{(x-y)^2}{(f_H(x)-f_H(y))^2}$. Note that $f_H\circ\varphi=f_{f_H(x),f_H(y),H}$. Therefore, taking $x,y$ to be $f_H(x),f_H(y)$ respectively in~\eqref{eq-restriction-SLE-arc-3}, we obtain~\eqref{eq-restriction-SLE-arc-2}.

We now integrate both sides of~\eqref{eq-restriction-SLE-arc-2} with respect to $\1_{\{0<x<1<y\}}\frac{dxdy}{|x-y|^2}$. The left side of~\eqref{eq-restriction-SLE-arc-2} becomes $\nu_{\hH\setminus H}(d\eta)$, while the right side becomes the same as the right side of~\eqref{eq-restriction-SLE-arc-1}. Indeed, if we perform a change of variables $u=f_H(x),v=f_H(y)$, then $\1_{\{0<x<1<y\}}\frac{|f'_H(x)f'_H(y)|}{|f_H(x)-f_H(y)|^2}dxdy=\1_{\{0<u<1<v\}}\frac{dudv}{(u-v)^2}$ as desired. This proves~\eqref{eq-restriction-SLE-arc-1}.

Finally, we multiply both sides of~\eqref{eq-restriction-SLE-arc-1} by $\frac{1}{Z}|g'_\eta(1)|$ (where $Z$ is defined in Lemma~\ref{lem-finite-SLE-arc}), then~\eqref{eq-restriction-SLE-arc} follows from $\frac{1}{Z}|g'_\eta(1)|\nu_{\hH\setminus H}(d\eta)=|f_H'(1)|\wt{\nu}_{\hH\setminus H}(d\eta)$.
\end{proof}

\subsection{Uniqueness of chordal hull}

In this Section, we prove Theorem~\ref{thm-chordal-hull} and Theorem~\ref{thm-unique-chordal-restriction}. Throughout this section, $J(z)=(\frac{1+z}{1-z})^2$ is a conformal map from $\D$ to $\C\setminus(-\infty,0]$ with $J(0)=1$.
We first consider the one-sided hull.
\begin{proof}[Proof of Theorem~\ref{thm-chordal-hull} (a)]
    Suppose $\PPP_1,\PPP_2$ are two probability measures on one-sided hull $H\in\mathcal{H}_+$ as in Theorem~\ref{thm-chordal-hull} (a). Let $F$ be a non-negative measurable function on $\eta$. We have
    \begin{equation*}
\int\left(\int \1_{\{H\cap\eta=\emptyset\}}e^{\frac{c}{2}m_\hH(H,\eta)}\PPP_1(dH)\right)F(\eta)\wt{\nu}_{\hH}(d\eta)=\int\left(\int \1_{\{H\cap\eta=\emptyset\}}e^{\frac{c}{2}m_\hH(H,\eta)}\PPP_2(dH)\right)F(\eta)\wt{\nu}_{\hH}(d\eta).
\end{equation*}

By Tonelli's theorem and Lemma~\ref{lem-restriction-SLE-arc}, we have
\begin{equation}
\label{eq-proof-one-side-1}
\int\left(\int F(\eta)\wt{\nu}_{\hH\setminus H}(d\eta)\right)|f_H'(1)|\PPP_1(dH)=\int\left(\int F(\eta)\wt{\nu}_{\hH\setminus H}(d\eta)\right)|f_H'(1)|\PPP_2(dH).
\end{equation}

Take $F=1$, we obtain from~\eqref{eq-proof-one-side-1} that $C:=\int|f_H'(1)|\PPP_1(dH)=\int|f_H'(1)|\PPP_2(dH)\in(0,1]$ (since $0<f'_H(1)\le1$). Let $\wt{\PPP}_1:=\frac{1}{C}|f_H'(1)|\PPP_1(dH)$ and $\wt{\PPP}_2:=\frac{1}{C}|f_H'(1)|\PPP_2(dH)$ be two probability measures. From~\eqref{eq-proof-one-side-1}, we obtain that
\begin{equation}
\label{eq-proof-one-side-2}
\iint F(f_H(\eta))\wt{\PPP}_1(dH)\wt{\nu}_\hH(d\eta)=\iint F(f_H(\eta))\wt{\PPP}_2(dH)\wt{\nu}_\hH(d\eta).
\end{equation}

Note that $f_H,g_\eta$ can be extended continuously to $\R_+$ and $f_H(\R_+)\subset\R_+,g_\eta(\R_+)\subset\R_+$. By Schwarz reflection theorem, let $\wh{f}_H$ (resp. $\wh{g}_\eta$) be the univalent map on $\C\setminus(-\infty,0]$ with $\wh{f}_H=f_H$ (resp. $\wh{g}_\eta=g_\eta$) on $\hH$. Note that $J^{-1}\circ\wh{f}_H\circ J,J^{-1}\circ\wh{g}_\eta\circ J\in\mathcal{D}$ (recall the space $\mathcal{D}$ in Section~\ref{section-composition}). Furthermore, the pushforward of $\wt{\PPP}_1$ (resp. $\wt{\PPP}_2$) via $H\mapsto J^{-1}\circ\wh{f}_H\circ J$ defines a probability measure $\mathsf{P}_1$ (resp. $\mathsf{P}_2$) on $\mathcal{D}$, and hence on $\mathcal{U}_1$.
The pushforward of $\wt{\nu}_\hH$ via $\eta\mapsto J^{-1}\circ\wh{g}_\eta\circ J$ defines another probability measure $\mathsf{Q}$. By~\eqref{eq-proof-one-side-2} and the arbitrariness of $F$, $(J^{-1}\circ\wh{f}_H\circ J)\circ(J^{-1}\circ\wh{g}_\eta\circ J)$ has the same law under $\mathsf{P}_1\times\mathsf{Q}$ or $\mathsf{P}_2\times\mathsf{Q}$. By Proposition~\ref{prop-unique-composition} (take $M=1$ there), we conclude that $\mathsf{P}_1=\mathsf{P}_2$. 
Finally, by Lemma~\ref{lem-measurable-chordal}, we have $\wt{\PPP}_1=\wt{\PPP}_2$ and hence $\PPP_1=\PPP_2$.
\end{proof}

We now consider the two-sided hull. 
To make notations clear, we write $\wt\nu_\hH^+:=\wt\nu_\hH$.
For a two-sided hull $H\in\mathcal{H}$, let $(\hH\setminus H)_+$ (resp. $(\hH\setminus H)_-$) be the connected component of $\hH\setminus H$ whose boundary contains $\R_+$ (resp. $\R_-$). Let $f_H^+$ (resp. $f_H^-$) be the unique conformal map from $\hH\to(\hH\setminus H)_+$ (resp. $\hH\to(\hH\setminus H)_-$) with $0,1,\infty$ (resp. $0,-1,\infty$) fixed.
To test the $(\hH\setminus H)_-$, we need to define the SLE arc measure on $\R_-$. Let $R(z):=-\overline{z}$ be the reflection map over the imaginary line $i\R$ and $\wt\nu_{\hH}^{-}:=R_{\#}\wt\nu_\hH^+$. For $\eta^-$ sampled from $\wt{\nu}_\hH^-$ with starting point $\eta_x^-$ and ending point $\eta_y^-$, let $g_\eta^-$ be the conformal map from $\hH$ to $D_{\eta^-}$ with $-1$ fixed and $g_\eta^-(0)=\eta_x^-$ and $g_\eta^-(\infty)=\eta_y^-$, where $D_{\eta^-}$ is the unique bounded connected component of $\hH\setminus\eta^-$.
By adapting the same proof as in Lemma~\ref{lem-restriction-SLE-arc}, the measure $\wt{\nu}_\hH^-$ satisfies the following conformal restriction identity, where $H\in\mathcal{H}$ and $\wt{\nu}_{\hH\setminus H}^-:=(f_H^-)_{\#}\wt{\nu}_\hH^-$:
\begin{equation}
\label{eq-restriction-SLE-arc-reflect}
|{f_H^{-}}'(-1)|\wt{\nu}_{\hH\setminus H}^-(d\eta^-)=\1_{\{\eta^-\cap H=\emptyset\}}e^{\frac{c}{2}m_\hH(H,\eta^-)}\wt{\nu}_\hH^-(d\eta^-).
\end{equation}

We can now prove the uniqueness of two-sided hull.
The main difference between the proof of one-sided and two-sided case is that we need to weight $\wt{\nu}_\hH^+(d\eta^+)\wt{\nu}_\hH^-(d\eta^-)$ by $\1_{\{\eta^+\cap\eta^-=\emptyset\}}e^{\frac{c}{2}m_\hH(\eta^+,\eta^-)}$.

\begin{proof}[Proof of Theorem~\ref{thm-chordal-hull} (b)]
    Suppose $\PPP_1,\PPP_2$ are two probability measures on two-sided hull $H\in\mathcal{H}$ as in Theorem~\ref{thm-chordal-hull} (b). Note that for $\eta^+$ (resp. $\eta^-$) sampled from $\wt\nu_\hH^+$ (resp. $\wt\nu_\hH^-$), if $\eta^+\cap\eta^-=\emptyset$, then $\overline{D_{\eta^+}\cup D_{\eta^-}}\in\mathcal{Q}$. Let $F$ be a non-negative measurable function on the pair $(\eta^+,\eta^-)$.
    We have
    \begin{equation*}
    \begin{aligned}
&\iint\left(\int\1_{\{H\cap(\eta^+\cup\eta^-)=\emptyset\}}e^{\frac{c}{2}m_\hH(H,\eta^+\cup\eta^-)}\PPP_1(dH)\right)F(\eta^+,\eta^-)\1_{\{\eta^+\cap\eta^-=\emptyset\}}e^{\frac{c}{2}m_\hH(\eta^+,\eta^-)}\wt{\nu}_\hH^+(d\eta^+)\wt{\nu}_\hH^-(d\eta^-)\\
&=\iint\left(\int\1_{\{H\cap(\eta^+\cup\eta^-)=\emptyset\}}e^{\frac{c}{2}m_\hH(H,\eta^+\cup\eta^-)}\PPP_2(dH)\right)F(\eta^+,\eta^-)\1_{\{\eta^+\cap\eta^-=\emptyset\}}e^{\frac{c}{2}m_\hH(\eta^+,\eta^-)}\wt{\nu}_\hH^+(d\eta^+)\wt{\nu}_\hH^-(d\eta^-).
\end{aligned}
\end{equation*}

Note that on the event $\eta^+\cap\eta^-=\emptyset$ and $H\cap(\eta^+\cup\eta^-)=\emptyset$, we have $m_\hH(H,\eta^+\cup\eta^-)+m_\hH(\eta^+,\eta^-)=m_\hH(H,\eta^+)+m_\hH(H,\eta^-)$. By Tonelli's theorem,~\eqref{eq-restriction-SLE-arc} and~\eqref{eq-restriction-SLE-arc-reflect}, we have
 \begin{equation*}
    \begin{aligned}
&\int \left(\iint F(\eta^+,\eta^-)\wt{\nu}^+_{\hH\setminus H}(d\eta^+)\wt{\nu}_{\hH\setminus H}^-(d\eta^-)\right)|{f_H^{+}}'(1)||{f_H^{-}}'(-1)|\PPP_1(dH)\\
&=\int \left(\iint F(\eta^+,\eta^-)\wt{\nu}^+_{\hH\setminus H}(d\eta^+)\wt{\nu}_{\hH\setminus H}^-(d\eta^-)\right)|{f_H^{+}}'(1)||{f_H^{-}}'(-1)|\PPP_2(dH).
\end{aligned}
\end{equation*}

Take $F=1$, we obtain that $C':=\int|{f_H^{+}}'(1)||{f_H^{-}}'(-1)|\PPP_1(dH)=\int|{f_H^{+}}'(1)||{f_H^{-}}'(-1)|\PPP_2(dH)\in(0,1]$ (since ${f_H^{+}}'(1),{f_H^{-}}'(-1)\in(0,1]$). Let $\wt\PPP_1(dH):=\frac{1}{C'}|{f_H^{+}}'(1)||{f_H^{-}}'(-1)|\PPP_1(dH)$ and $\wt\PPP_2(dH):=\frac{1}{C'}|{f_H^{+}}'(1)||{f_H^{-}}'(-1)|\PPP_2(dH)$ be two probability measures. We have
\begin{equation}
\label{eq-proof-chordal-hull-1}
\begin{aligned}
&\iiint F(f_H^+(\eta^+),f_H^-(\eta^-))\wt\PPP_1(dH)\wt{\nu}^+_{\hH}(d\eta^+)\wt{\nu}_{\hH}^-(d\eta^-)\\
&=\iiint F(f_H^+(\eta^+),f_H^-(\eta^-))\wt\PPP_2(dH)\wt{\nu}^+_{\hH}(d\eta^+)\wt{\nu}_{\hH}^-(d\eta^-).
\end{aligned}
\end{equation}

Recall the reflection map $R(z)=-\overline{z}$. By Schwarz reflection theorem, let $\wh{f}_H^+,\wh{g}_\eta^+$ (resp. $\wh{f}_H^-,\wh{g}_\eta^-$) be the univalent map on $\C\setminus(-\infty,0]$ (resp. $\C\setminus[0,+\infty)$) such that $\wh{f}_H^+=f_H^+,\wh{g}_\eta^+=g_\eta^+$ (resp. $\wh{f}_H^-=f_H^-,\wh{g}_\eta^-=g_\eta^-$) on $\hH$. The pushforward of $\wt\PPP_1$ (resp. $\wt\PPP_2$) via $H\mapsto (J^{-1}\circ \wh{f}_H^+\circ J,J^{-1}\circ R^{-1}\circ \wh{f}_H^-\circ R\circ J)$ defines a probability measure $\mathsf{P}_1$ (resp. $\mathsf{P}_2$) on $\mathcal{D}\times\mathcal{D}$, and hence on $\mathcal{U}_1\times\mathcal{U}_1$. Furthermore, the pushforward of $\wt{\nu}^+_{\hH}(d\eta^+)\times\wt{\nu}_{\hH}^-(d\eta^-)$ via $(\eta^+,\eta^-)\mapsto (J^{-1}\circ \wh{g}_\eta^+\circ J,J^{-1}\circ R^{-1}\circ \wh{g}_\eta^-\circ R\circ J)$ gives a probability measure $\mathsf{Q}\times\mathsf{Q}$ on $\mathcal{D}\times\mathcal{D}$. Now,~\eqref{eq-proof-chordal-hull-1} and the arbitrariness of $F$ implies that
\begin{equation*}
\big((J^{-1}\circ \wh{f}_H^+\circ J)\circ(J^{-1}\circ \wh{g}_\eta^+\circ J),(J^{-1}\circ R^{-1}\circ \wh{f}_H^-\circ R\circ J)\circ(J^{-1}\circ R^{-1}\circ \wh{g}_\eta^-\circ R\circ J)\big)
\end{equation*}
has the same law under $\mathsf{P}_1\times\mathsf{Q}\times\mathsf{Q}$ or $\mathsf{P}_2\times\mathsf{Q}\times\mathsf{Q}$. By the $n=2$ case of Proposition~\ref{prop-joint-composition}, we conclude that $\mathsf{P}_1=\mathsf{P}_2$.
Finally, by Lemma~\ref{lem-measurable-chordal}, we have $\wt{\PPP}_1=\wt{\PPP}_2$ and hence $\PPP_1=\PPP_2$.
\end{proof}

We now show the uniqueness of chordal generalized restriction measure for $c\le1$.
\begin{proof}[Proof of Theorem~\ref{thm-unique-chordal-restriction}]
    From Definition~\ref{def-chordal-restriction}, for a one-sided (resp. two-sided) restriction measure with exponent $\alpha$ and $A\in\mathcal{Q}_+$ (resp. $A\in\mathcal{Q}$), we have $\E_\PPP[\1_{\{H\cap A=\emptyset\}}e^{\frac{c}{2}m_\hH(H,A)}]=(\phi_A'(0)\phi_A'(\infty))^\alpha$. Therefore, Theorem~\ref{thm-chordal-hull} implies the uniqueness of such $\PPP$. For the existence part, the $\SLE_\kappa(\rho)$ with $\rho>-2$ from 0 to $\infty$ with force point $0_-$ gives the one-sided restriction measure with exponent $\alpha=\frac{(\rho+2)(\rho+6-\kappa)}{4\kappa}$ for all $\alpha>0$; see~\cite{lsw-restriction} and also~\cite[Section 6.2.1]{qian2021generalized}. Note that the $\alpha=0$ case was given by the deterministic hull $(-\infty,0]$.
    The construction of two-sided restriction measure with exponent $\alpha\ge\frac{6-\kappa}{2\kappa}$ was also given in~\cite[Proposition 6.2]{qian2021generalized}.

    It remains to show the non-existence part.  First, we must have $\alpha\ge0$ since $\phi_A'(0)\phi_A'(\infty)\in(0,1]$ for $A\neq\emptyset$. Indeed, take $F=1$ in~\eqref{eq-proof-one-side-1}, we obtain that $\int (\phi_{D_\eta}'(0)\phi_{D_\eta}'(\infty))^\alpha\wt\nu_\hH(d\eta)=\int |f_H'(1)|\PPP(dH)\in(0,1]$ since $f_H'(1)\in(0,1]$. 
    Second, if $\alpha=0$, then by the one-sided case, the right (resp. left) boundary of the hull is $(-\infty,0]$ (resp. $[0,+\infty$)), a contradiction.
    It remains to rule out the existence of two-sided restriction measure with exponent $\alpha\in(0,\frac{6-\kappa}{2\kappa})$.
    The following method is the same as~\cite[Corollary 8.6]{lsw-restriction} and we will keep brief. For $H$ sampled from such two-sided restriction measure $\PPP$ with exponent $\alpha\in[0,\frac{6-\kappa}{2\kappa})$, the right boundary of $H$ is an $\SLE_\kappa(\rho)$ with force point $0_-$ and $\rho<0$. Therefore, $\PPP[i\in(\hH\setminus H)_+]>\frac{1}{2}$. Similarly, we have $\PPP[i\in(\hH\setminus H)_-]>\frac{1}{2}$, which leads to a contradiction.
\end{proof}

\begin{remark}
    \label{remark-trichordal}
    Our approach can also be applied to the uniqueness of generalized tri-chordal restriction measure with central charge $c\le1$; see~\cite{Qian-trichordal} and~\cite[Section 6.2.2]{qian2021generalized} for more background. We only outline the necessary changes in proof and leave the details to interesting readers. For a trichordal hull $H$ on $\hH$ with $H\cap\partial\hH=\{0,1,\infty\}$, let $(\hH\setminus H)_a,(\hH\setminus H)_b,(\hH\setminus H)_c$ be the three connected components of $\hH\setminus H$ whose boundaries contain $(0,1),(1,\infty),(-\infty,0)$ respectively. We may define $\wt\nu_\hH^a,\wt\nu_\hH^b,\wt\nu_\hH^c$ to be the SLE arc probability measures on $(0,1),(1,\infty),(-\infty,0)$ respectively. 
    We use the following measure to test against $H$:
    \begin{equation*}
\1_{\{\eta^a\cap\eta^b=\emptyset,\eta^b\cap\eta^c=\emptyset,\eta^c\cap\eta^a=\emptyset\}}
e^{\frac{c}{2}(m_\hH(\eta^a,\eta^b)+m_\hH(\eta^a\cup\eta^b,\eta^c))}\wt\nu_\hH^a(d\eta^a)\wt\nu_\hH^b(d\eta^b)\wt\nu_\hH^c(d\eta^c).
\end{equation*}
On the event $H\cap(\eta^a\cup\eta^b\cup\eta^c)=\emptyset$,
we have
\begin{equation*}
m_\hH(H,\eta^a\cup\eta^b\cup\eta^c)+m_\hH(\eta^a,\eta^b)+m_\hH(\eta^a\cup\eta^b,\eta^c)=\sum_{i\in\{a,b,c\}}m_\hH(H,\eta^i).
\end{equation*}
We can then apply the same strategy, and use the $n=3$ case of Proposition~\ref{prop-joint-composition} to complete the proof.
\end{remark}

\subsection{Uniqueness of radial hull}
\label{subsection-radial}
In this Section, we prove the uniqueness result for radial hull. It is convenient for us to work on the unit disk $\D$ with marked points $0,1$. The proof is almost the same as in the chordal one-sided case, except that the topology is different. 

Let $\mathcal{H}_\D$ be the collection of compact sets $H\subset\overline{\D}$ such that $0\in H,H\cap\partial\D=\{1\}$ and $\D\setminus H$ is simply connected. We call an $H\in\mathcal{H}_\D$ a radial hull. Let $\mathcal{Q}_\D$ be the collection of compact sets $A\subset\overline{\D}$ with $0,1\notin A$ and $\D\setminus A$ is simply connected. We equip $\mathcal{H}_\D$ with $\sigma$-algebra generated by $\{H:H\cap A=\emptyset\}$ for all $A\in\mathcal{Q}_\D$. For two subsets $A_1,A_2\subset\overline{\D}$, let $m_{\D}(A_1,A_2)$ be the total mass of the Brownian loops on $\D$ that intersect both $A_1$ and $A_2$.
\begin{theorem}
    \label{thm-radial-hull}
    Fix $c\le 1$. Suppose $\PPP_1,\PPP_2$ are two probability measures on $\mathcal{H}_\D$ satisfying
    \begin{equation*}
\E_{\PPP_1}[\1_{\{H\cap A=\emptyset\}}e^{\frac{c}{2}m_\D(H,A)}]=\E_{\PPP_2}[\1_{\{H\cap A=\emptyset\}}e^{\frac{c}{2}m_\D(H,A)}]\in[0,\infty]
\end{equation*}
    for every $A\in\mathcal{Q}_\D$. Then $\PPP_1=\PPP_2$.
\end{theorem}

To show Theorem~\ref{thm-radial-hull}, we first introduce
the disk analogue of the SLE arc measure $\nu_\hH$. For $0<u<\pi<v<2\pi$, let $\SLE_{\D;e^{iu}\to e^{iv}}$ be the $\SLE_\kappa(\kappa-2)$ probability measure on $\D$ from $e^{iu}$ to $e^{iv}$ with force point $e^{iu_+}$. We define
\begin{equation*}
\nu_\D(d\eta):=
\iint_{0<u<\pi<v<2\pi}
\SLE_{\D;e^{iu}\to e^{iv}}(d\eta)
\frac{du dv}{|e^{iu}-e^{iv}|^2}.
\end{equation*}

Here, we restrict on $0<u<\pi<v<2\pi$ since we choose $-1=e^{i\pi}$ as a marked point for $\nu_\D$ (which is similar to that the marked point of $\nu_\hH$ is 1).
The kernel $\frac{du dv}{|e^{iu}-e^{iv}|^2}$ is just the boundary Poisson kernel on $\partial\D$.
For $\nu_\D$-a.e. $\eta$, let $A_\eta$ be the connected component of $\D\setminus\eta$ whose boundary contains $-1$. Let $h_\eta:\hH\to A_\eta$ be conformal with $h_\eta(0)=e^{iu},h_\eta(1)=-1,h_\eta(\infty)=e^{iv}$. Recall the constant $Z$ in Lemma~\ref{lem-finite-SLE-arc}. We define a measure (not a probability measure)
\begin{equation*}
\wt\nu_\D(d\eta):=
\frac{1}{Z}\1_{\{0\notin \overline{A_\eta}\}}|h_\eta'(1)|\nu_\D(d\eta).
\end{equation*}
For $H\in\mathcal{H}_\D$, let $f_H:\hH\to\D\setminus H$ be conformal with $f_H(0)=1_+,f_H(1)=-1,f_H(\infty)=1_-$.

\begin{lemma}
    \label{lem-restriction-SLE-arc-radial}
    For every $H\in\mathcal{H}_\D$, we have
    \begin{equation}
    \label{eq-restriction-SLE-arc-radial}
|f_H'(1)|(f_H)_\#\wt\nu_\hH(d\eta)
=
\1_{\{H\cap\eta=\emptyset\}}
e^{\frac{c}{2}m_\D(H,\eta)}
\wt\nu_\D(d\eta).
    \end{equation}
\end{lemma}
\begin{proof}
Similar to Lemma~\ref{lem-restriction-SLE-arc} and the proof of~\eqref{eq-restriction-SLE-arc-1}, we have
\begin{equation}
\label{eq-restriction-SLE-arc-radial-1}
(f_H)_\#\nu_\hH(d\eta)=\1_{\{\eta\cap H=\emptyset\}}e^{\frac{c}{2}m_\D(H,\eta)}\nu_\D(d\eta).
\end{equation}
Note that $\{\eta\cap H=\emptyset\}$ implies $\{0\notin A_\eta\}$. Therefore,~\eqref{eq-restriction-SLE-arc-radial} follows by multiplying both sides of~\eqref{eq-restriction-SLE-arc-radial-1} by $\frac{1}{Z}|h_\eta'(1)|$, as well as the factorization $|h_\eta'(1)|=|f'_H(1)|\cdot|g'_{f_H^{-1}(\eta)}(1)|$ and $\frac{1}{Z}|g'_{f_H^{-1}(\eta)}(1)|(f_H)_\#\nu_\hH(d\eta)=(f_H)_\#\wt\nu_\hH(d\eta)$.
\end{proof}

We can now prove the uniqueness of radial hull. Recall $J(z)=(\frac{1+z}{1-z})^2$ is conformal from $\D$ to $\C\setminus (-\infty,0]$ with $J(0)=1$.

\begin{proof}[Proof of Theorem~\ref{thm-radial-hull}]

Note that $\overline{A_\eta}\in\mathcal{Q}_\D$ for $\wt\nu_\D$-a.e. $\eta$. For a non-negative measurable function $F(\eta)$, we have
\begin{equation*}
\int\left(\int\1_{\{H\cap \eta=\emptyset\}}e^{\frac{c}{2}m_\D(H,\eta)}\PPP_1(dH)\right)F(\eta)\wt\nu_\D(d\eta)=\int\left(\int\1_{\{H\cap \eta=\emptyset\}}e^{\frac{c}{2}m_\D(H,\eta)}\PPP_2(dH)\right)F(\eta)\wt\nu_\D(d\eta).
\end{equation*}

By Tonelli's theorem and Lemma~\ref{lem-restriction-SLE-arc-radial}, we have
\begin{equation}
\label{eq-proof-radial-1}
\iint F(f_H(\eta))|f_H'(1)|\PPP_1(dH)\wt\nu_\hH(d\eta)=
\iint F(f_H(\eta))|f_H'(1)|\PPP_2(dH)\wt\nu_\hH(d\eta).
\end{equation}

By Schwarz reflection theorem, $f_H$ extends to a univalent map (since $f_H(z)\neq0$) on $\C\setminus(-\infty,0]$ given by $\wh{f}_H(z):=f_H(z)$ if $\Im(z)>0$ and $\wh{f}_H(z):=1/\overline{f_H(\overline{z})}$ if $\Im(z)<0$.
We define $F_H(z):=i(\wh f_H(J(z))+1)$ to be univalent from $\D$ to $\C$. Note that $F_H(0)=0$ and $F_H'(0)>0$ (since $f'_H(1)=-i|f_H'(1)|$).
Moreover, since $i\notin F_H(\D)$, Koebe's quarter theorem implies that $4|f_H'(1)|=F_H'(0)\le4$, and hence $|f_H'(1)|\in(0,1],F_H\in\mathcal{U}_4$. (recall the space $\mathcal{U}_M$ in Section~\ref{section-composition}).

Take $F=1$, we obtain that $C'':=\int |f_H'(1)|\PPP_1(dH)=\int|f_H'(1)|\PPP_2(dH)\in(0,1]$.
Let $\wt\PPP_i(dH):=\frac{1}{C''}|f_H'(1)|\PPP_i(dH)$ for $i=1,2$. By~\eqref{eq-proof-radial-1},
\begin{equation}
\label{eq-proof-radial-2}
\iint F(f_H(\eta))\wt\PPP_1(dH)\wt\nu_\hH(d\eta)
=
\iint F(f_H(\eta))\wt\PPP_2(dH)\wt\nu_\hH(d\eta).
\end{equation}

Recall the reflected map $\wh g_\eta$ in the proof of Theorem~\ref{thm-chordal-hull} (a). For $i=1,2$, the pushforward of $\wt\PPP_i$ under $H\mapsto F_H$ defines a probability measure $\mathsf{P}_i$ on $\mathcal{U}_4$.
The pushforward of $\wt\nu_\hH$ under $\eta\mapsto J^{-1}\circ\wh g_\eta\circ J$ defines a probability measure $\mathsf{Q}$ on $\mathcal{D}$. Since $F_H(J^{-1}(\wh{g}_\eta(J(z))))=i(\wh{f}_H(\wh{g}_\eta(J(z)))+1)$,~\eqref{eq-proof-radial-2} implies that $F_H\circ( J^{-1}\circ\wh g_\eta\circ J)$ has the same law under $\mathsf{P}_{1}\times\mathsf{Q}$ or $\mathsf{P}_{2}\times\mathsf{Q}$. 
By Proposition~\ref{prop-unique-composition}, we conclude that $\mathsf{P}_1=\mathsf{P}_2$.
By Lemma~\ref{lem-measurable-radial}, we have $\wt\PPP_1=\wt\PPP_2$ and hence $\PPP_1=\PPP_2$.
\end{proof}

A direct consequence of Theorem~\ref{thm-radial-hull} is the uniqueness of generalized radial restriction measure which we now recall; see~\cite{Wu-radial-restriction} for the study of $c=0$ and~\cite{qian2021generalized} for $c\le1$.
\begin{definition}[{\cite[Definition 1.3]{qian2021generalized}}]
    \label{def-radial-restriction}
    Suppose $\kappa\in(0,4],c\le1$ are related by~\eqref{eq-c-kappa}. We say a probability measure $\PPP$ on $\mathcal{H}_\D$ satisfies generalized radial restriction property with central charge $c$ and exponent $(\alpha,\beta)$, if
    \begin{equation*}
    \1_{\{H\cap A=\emptyset\}}\PPP(dH)=|\phi_A'(0)|^\alpha|\phi_A'(1)|^\beta e^{-\frac{c}{2}m_{\D}(H,A)}\PPP_A(dH)
\end{equation*}
for all $A\in\mathcal{Q}_\D$. Here $\phi_A:\D\setminus A\to\D$ is conformal with $0,1$ fixed and $\PPP_A:=(\phi_A^{-1})_\#\PPP$.
\end{definition}

In particular, for a generalized radial restriction measure $\PPP$ with central charge $c$ and exponent $(\alpha,\beta)$, we have $\E_\PPP[\1_{\{H\cap A=\emptyset\}}e^{\frac{c}{2}m_\D(H,A)}]=|\phi_A'(0)|^\alpha|\phi_A'(1)|^\beta$ for all $A\in\mathcal{Q}_\D$. Therefore, Theorem~\ref{thm-radial-hull} leads to the following result.
\begin{theorem}
    \label{thm-uniqueness-radial-restriction}
    Fix $c\le 1$ and $\kappa\in(0,4]$ related to $c$ via~\eqref{eq-c-kappa}.
    For each $(\alpha,\beta)$, there is at most one generalized radial restriction measure with central charge $c$ and exponent $(\alpha,\beta)$.
\end{theorem}

In~\cite[Theorem 1.6]{qian2021generalized}, Qian used radial hypergeometric SLE to construct the generalized radial restriction measure for $\alpha\le\eta_\kappa(\beta)$ and $\beta\ge\frac{6-\kappa}{2\kappa}$, where $\eta_\kappa(\beta):=\frac{(\sqrt{16\beta\kappa+(4-\kappa)^2}-(4-\kappa))^2-4(4-\kappa)^2}{32\kappa}$ is the generalized disconnection exponent. 
The optimality of this range for $(\alpha,\beta)$ was conjectured in~\cite{qian2021generalized} and proved only for $c=0$ in~\cite{Wu-radial-restriction}. The following observation supports this conjecture.
Heuristically, zooming in near $1$ should yield a generalized two-sided chordal restriction measure with exponent $\beta$. Theorem~\ref{thm-unique-chordal-restriction} then forces that $\beta\ge\frac{6-\kappa}{2\kappa}$. Moreover, the $\alpha=\eta_\kappa(\beta)$ is critical in the sense that $0\in\partial H$ almost surely; see~\cite[Theorem 1.6 (iv)]{qian2021generalized}. Since our main curiosity lies in the uniqueness result, we do not address this conjecture here.
\begin{remark}
    \label{remark-weight-conformal-radius}
    We explain how weighting by a power of the conformal radius changes the exponent $\alpha$. Fix $c\le1$. Let $\PPP^{(\alpha,\beta)}$ be the generalized radial restriction measure with exponent $(\alpha,\beta)$.
    Suppose $\alpha<\eta_\kappa(\beta)$ and $\beta\ge\frac{6-\kappa}{2\kappa}$. As demonstrated in~\cite[Theorem 1.6]{qian2021generalized}, in this case $0$ lies in the interior of $H$ for a.s. $H$ sampled from $\PPP^{(\alpha,\beta)}$. Let $D_0$ be the connected component of $\D\setminus\partial H$ that contains 0. Let $\varphi_H:\D\to D_0$ be conformal with $\varphi_H(0)=0$. $|\varphi_H'(0)|\in(0,1)$ is also known as the conformal radius of $D_0$ seen from 0.
    For $t\in\R$, let $Z_t^{(\alpha,\beta)}:=\int|\varphi_H'(0)|^t\PPP^{(\alpha,\beta)}(dH)\in(0,\infty]$. It is direct to verify that $\frac{1}{Z_t^{(\alpha,\beta)}}|\varphi_H'(0)|^t\PPP^{(\alpha,\beta)}(dH)$ satisfies generalized radial restriction measure with exponent $(\alpha-t,\beta)$ if $Z_t^{(\alpha,\beta)}<\infty$. Therefore, Theorem~\ref{thm-uniqueness-radial-restriction} yields the following:
    \begin{equation*}
\frac{d\PPP^{(\alpha-t,\beta)}(H)}{d\PPP^{(\alpha,\beta)}(H)}=\frac{1}{Z_t^{(\alpha,\beta)}}|\varphi_H'(0)|^t.
\end{equation*}
Also, one can see from~\cite[Theorem 1.5]{qian2021generalized} that $Z_t^{(\alpha,\beta)}<\infty$ if and only if $t>\alpha-\eta_\kappa(\beta)$. 
\end{remark}

\appendix
\section{Measurability}
\label{appendix-measurable}
For $H\in\mathcal{H_+}$, recall that $f_H:\hH\to \hH\setminus H$ is conformal with $0,1,\infty$ fixed, and $\wh f_H$ be its Schwarz reflection along $\R_+$. Recall $J:\D\to\C\setminus(-\infty,0]$ is defined as $J(z)=(\frac{1+z}{1-z})^2$. 

\begin{lemma}
\label{lem-measurability}
The map $T:\mathcal{H}_+\to\mathcal{D}$ defined by
$H\mapsto J^{-1}\circ\wh f_H\circ J$ is a measurable map.
\end{lemma}
\begin{proof}
Recall that we equip $\mathcal{D}$ with the Borel $\sigma$-algebra generated by the topology of locally uniform convergence. Recall that the $\sigma$-algebra of $\mathcal{H}_+$ is generated by $\{H:H\cap A=\emptyset\}$ for all $A\in\mathcal{Q}_+$. 
We use the following criterion to show the measurability of $T$: it suffices to find a sequence of measurable maps $T_n:\mathcal{H}_+\to\mathcal{D}$, such that for each $H\in\mathcal{H_+}$, $T_n(H)\to T(H)$ as $n\to\infty$.

For a set $K\subset\C$, let
$K^*:=\{\overline z:z\in K\}$.
From now on, we fix an $H\in\mathcal{H}_+$ and explain how to define $T_n$.
Let $(A_j)_{j\ge1}$ be the collection of all the closed polygonal
hulls in $\mathcal{Q}_+$ with rational vertices and
$A_j\cap\R=[x_j,y_j]$, where $0<x_j<1<y_j$.
We have $\C\setminus(H\cup H^*)=
\cup_{j:H\cap A_j=\emptyset}(A_j\cup A_j^*)^\circ$, where $(A_j\cup A_j^*)^\circ$ is the interior of $A_j\cup A_j^*$.
For each $n\ge1$, let $V_n(H):=
\cup_{{j\le n:H\cap A_j=\emptyset}}(A_j\cup A_j^*)^\circ$.
When $V_n(H)$ is nonempty, let $\mathcal{F}(V_n(H))$ be its filling
obtained by adding all bounded complementary components.
Note that $\mathcal{F}(V_n(H))\subset\C\setminus(H\cup H^*)$ and $\mathcal{F}(V_n(H))$ eventually increases to $\C\setminus(H\cup H^*)$.

Let $\psi_{n,H}:\D\to\mathcal{F}(V_n(H))$ be conformal with
$\psi_{n,H}(0)=1$ and $\psi_{n,H}'(0)>0$.
When $V_n(H)$ is empty, we define $\psi_{n,H}=J$.
For each $n\ge1$, let $T_n:\mathcal{H}_+\to\mathcal{D}$ be the map
$H\mapsto J^{-1}\circ\psi_{n,H}$.
Note that $T_n$ is measurable since it depends only on the
measurable events $\{H:H\cap A_j=\emptyset\}$ for $j\le n$.

By the Carath\'eodory kernel theorem,
$\psi_{n,H}\to\wh f_H\circ J$ locally uniformly on $\D$.
Consequently, $J^{-1}\circ\psi_{n,H}
\to J^{-1}\circ\wh f_H\circ J$
locally uniformly, and hence $T_n(H)\to T(H)$ as $n\to\infty$.
\end{proof}

By Lemma~\ref{lem-measurability}, for a probability measure $\PPP$ on $\mathcal{H_+}$, we can define the pushforward measure $\mathsf{P}:=T_\#\PPP$ on $\mathcal{D}$. Note that $T$ is a Borel injection from $\mathcal{H}_+$ to $\mathcal{D}$.
\begin{lemma}
    \label{lem-measurable-chordal}
     $\PPP$ is determined by $\mathsf{P}$.
\end{lemma}
\begin{proof}
For $A\in\mathcal{Q}_+$, consider the collection of conformal maps $\mathcal{D}_A:=\{f\in\mathcal{D}:J^{-1}(A\cup A^*)\subset f(\D)\}$.
By Schwarz reflection, we have $\{H:H\cap A=\emptyset\}=T^{-1}(\mathcal{D}_A)$.
Note that $\mathcal{D}_A$ is an open set under locally uniform convergence, since $J^{-1}(A\cup A^*)$ is compact.
In particular, we have $\mathscr{F}_{\mathcal{H}_+}\subset T^{-1}(\mathscr{F}_{\mathcal{D}})$, where $\mathscr{F}_{\mathcal{H}_+}$ (resp. $\mathscr{F}_{\mathcal{D}}$) is the $\sigma$-algebra of $\mathcal{H}_+$ (resp. $\mathcal{D}$). Hence 
$\PPP$ is determined by $\mathsf{P}$.
\end{proof}

We also need a radial version of Lemma~\ref{lem-measurability} and Lemma~\ref{lem-measurable-chordal}. 
Recall that in Section~\ref{subsection-radial}, for $H\in\mathcal{H}_\D$, we define $f_H:\hH\to\D\setminus H$ as a conformal map with $f_H(0)=1_+,f_H(1)=-1,f_H(\infty)=1_-$. Recall the reflection map $\wh f_H$ and $F_H:=i(\wh f_H(J(z))+1)\in\mathcal{U}_4$ defined below~\eqref{eq-proof-radial-1}.
\begin{lemma}
    \label{lem-measurable-radial}
     The map $T_\D:\mathcal{H}_\D\to\mathcal{U}_4$ defined by
$H\mapsto F_H$ is a measurable and injective map. Furthermore, each probability law $\PPP$ on $\mathcal{H}_\D$ is determined by its pushforward law $(T_\D)_\#\PPP$.
\end{lemma}

We omit its proof since its identical to that of Lemma~\ref{lem-measurability} and Lemma~\ref{lem-measurable-chordal}. 

\small
\def\cprime{$'$}


\end{document}